\documentclass[11pt]{article}
\usepackage{mathtools}
\usepackage{amsthm, amssymb}
\usepackage{setspace, bigints}
\usepackage{mathrsfs}
\usepackage{mleftright}
\usepackage{bm}
\usepackage{bbm}
\usepackage[margin=0.75 in]{geometry}
\usepackage{caption}
\usepackage{subcaption}
\usepackage[toc,page]{appendix}     
\usepackage{pdfpages}
\usepackage{epstopdf}
\usepackage{booktabs}
\usepackage{graphicx}
\usepackage{listings}
\usepackage{xcolor}
\usepackage{tabularx, multirow}
\usepackage[ruled,linesnumbered]{algorithm2e}
\usepackage{tabu}
\usepackage{comment}
\usepackage{natbib}
\usepackage{authblk}

\usepackage{amsthm}

\newtheorem{assumptiontemp}{Assumption}

\newenvironment{Assumption}
  {\begin{assumptiontemp}}
  {\end{assumptiontemp}}

\usepackage{empheq}
\usepackage{accents}
\usepackage{bm}
\newcommand*{\B}[1]{\ifmmode\bm{#1}\else\textbf{#1}\fi}

\newcommand{\be} {\begin{eqnarray*}}
\newcommand{\ee} {\end{eqnarray*}}

\DeclareMathOperator{\KL}{KL}
\DeclareMathOperator{\Law}{Law}

\newcommand{\argmin}{\mathop{\rm argmin~}}

\newcommand{\dd}{{\rm d}}

\newcommand{\wht}{\widehat}

\def\m{\mathcal}
\def\ms{\mathscr}
\def\mb{\mathbb}

\def\mx{\mbox}

\newcommand{\matnorm}[1]{{\left\vert\kern-0.25ex\left\vert\kern-0.25ex\left\vert #1 
    \right\vert\kern-0.25ex\right\vert\kern-0.25ex\right\vert}}

\newtheorem{theorem}{Theorem}[section]

\newtheorem{proposition}[theorem]{Proposition}

\newtheorem{rem}{Remark}

\numberwithin{equation}{section}

\usepackage{hyperref}
\usepackage[capitalize]{cleveref}

\usepackage{xspace}

\title{Convergence Rate of the Solution of Multi-marginal Schrodinger Bridge Problem with Marginal Constraints from SDEs}

\author{Rentian Yao}
\author{Young--Heon Kim}
\author{Geoffrey Schiebinger}
\affil{Department of Mathematics, University of British Columbia \authorcr Email: \{rentian2, yhkim, geoff\}@math.ubc.ca}
\date{\vspace{-2em}}
\date{}

\begin{document}
\maketitle
\begin{abstract}
In this paper, we investigate the multi-marginal Schrodinger bridge (MSB) problem whose marginal constraints are marginal distributions of a stochastic differential equation (SDE) with a constant diffusion coefficient, and with time dependent drift term. As the number $m$ of marginal constraints increases, we prove that the solution of the corresponding MSB problem converges to the law of the solution of the SDE at the rate of $O(m^{-1})$, in the sense of KL divergence.  
Our result extends the work of~\cite{agarwal2024iterated} to the case where the drift of the underlying stochastic process is time-dependent.
\end{abstract}
\section{Introduction}
In recent years, there has been increasing interest in Schrodinger bridge problems, with applications in trajectory inference~\citep{lavenant2024toward, chizat2022trajectory, yao2025learning, hong2025trajectory, gu2025private}, 
computational optimal transport~\citep{cuturi2013sinkhorn, li2025multimarginal, knight2008sinkhorn}, distribution matching~\cite{liu2023generalized, yang2025topological}, score estimation~\cite{agarwal2024iterated, mordant2024entropic}, and diffusion generative models~\citep{deng2024reflected, de2021diffusion, wang2021deep}.

Motivated by the finite-sample analysis of trajectory inference~\cite{yao2025learning}, this paper focuses on investigating the explicit convergence rate for the multi-marginal Schrodinger bridge problem as the number of marginal constraints $m$, equivalently the number of time sample points, increases, in the case where the marginals are the distributions from a stochastic differential equation (SDE) given at different time points; see Theorem \ref{thm: main result}.
Under this SDE assumption, we show that the KL divergence between the path-space distribution of this SDE and the multi-marginal Schrodinger bridge with $m$ marginal constraints converges at the rate of $O(m^{-1})$.
Our result extends the work of~\cite{agarwal2024iterated} to the case where the drift of the underlying stochastic process is time-dependent. This generalization is particularly important in practical applications such as developmental biology, where time-dependent drift coefficients commonly arise~\cite{farrell2018single, lavenant2024toward, schiebinger2019optimal}.

\subsection{Multi-marginal Schrodinger bridge}
Before we present the problem, let us first fix our notation:
\begin{itemize}
    \item Let $\m X$ be the state space, either $\m X= \mathbb{R}^d$, the flat torus $\mathbb{T}^d = [0, 2\pi)^d$, or a Riemannian manifold, and let $\lambda_{\m X}$ be the Lebesgue (volume) measure on $\m X$. 
\item For a fixed $T > 0$, define the path space $\Omega := \m C([0, T]; \m X)$ as the set of all continuous maps from $[0, T]$ to $\m X$. 
Also, let
\begin{align*}
    \hbox{$\Omega_{[s, t]} :=\m C([s, t]; \m X)$ for $0\le s< t \le T$}.
\end{align*}

\item For any stochastic process $\{X_t: t\in[0, T]\}$ with path-space distribution $R\in\ms P(\Omega)$, let 
$$R_{t_1, t_2} = \Law(X_{t_1}, X_{t_2})$$ 
denote the joint distribution of $X_{t_1}$ and $X_{t_2}$, 
$$R_{[t_1, t_2]} = \Law(\{X_t: t\in[t_1, t_2]\})$$
represent the law of the process restricted to the time interval $[t_1, t_2]$ for any $0\leq t_1 < t_2\leq T$, and 
$$R_{t_2\,|\,t_1} = \Law(X_{t_2}\,|\,X_{t_1})$$
be the conditional distribution of $X_{t_2}$ given $X_{t_1}$.
Notice that $R_{t_1, t_2} \in \ms P (\m X \times \m X)$, $R_{[t_1, t_2]} \in \ms P(\Omega_{[t_1, t_2]})$, and $R_{t_2\,|\,t_1}\in\ms P(\m X)$ for any given $X_{t_1}$.
For more time points, $t_1, ..., t_k$,  we similarly define $R_{t_1, ..., t_k} = {\rm Law} (X_{t_1}, ..., X_{t_k}).$

\item For a given $\tau > 0$, let $W^\tau\in\ms P(\Omega)$ be the law of the reversible Brownian motion on $\m X$ with temperature $\tau$.
\item Given two probability measures $P$ and $Q$ on a probability space, Kullback--Leibler (KL) divergence is given by
\begin{align*}
    \KL(P\,\|\, Q) = \mb E_{P}\Big[\log\frac{\dd P}{\dd Q} \Big]= \int  \log\left( \frac{\dd P}{\dd Q}\right) \dd P.
\end{align*}
(Note that by Jensen's inequality it holds that $ \KL(P\,\|\, Q) \ge 0.$)
\item We let $\ms P_2^r(\m X)$ denote the set of all probability measures $\mu$ on $\m X$ that are absolutely continuous with respect to $\lambda_{\m X}$, and with finite second moment, that is, $\int_{\m X} dist^2 (x_0, x)d\mu (x)$ for a fixed $x_0 \in \m X$.

\end{itemize}

For any $m+1$ probability distributions $\mu_0, \mu_1, \dots, \mu_m\in\ms P(\m X)$,  
and for a sequence of time points $\m T_m = (t_0, t_1, \dots, t_m)\in(0,T)^m$
, we define the multi-marginal Schrodinger bridge problem (MSB) as 
the following constrained entropy minimization problem,
\begin{align}\label{eqn: general_MSB}
\wht R^{\m T_m} = \argmin_{\substack{R\in\ms P(\Omega)\\ R_{t_j} = \mu_j, j=0,\dots, m}} \KL(R\,\|\,W^\tau).
\end{align} 
Notice that $R\mapsto \KL(R\, \|\, W^\tau)$  is a strictly convex function, thus \eqref{eqn: general_MSB} has a unique solution if exists.

As the KL-diverence gives a quantitative comparison between two probability distributions, the solution $\wht R^{\m T_m}$ to \eqref{eqn: general_MSB} gives a closest stochastic process to the Brownian motion $W^\tau$, among those satisfying the given marginal condition ($R_{t_j} = \mu_j$).  
This is relevant to many practical problems where one seeks to estimate a process from marginals, including trajectory inference~\cite{lavenant2024toward, yao2025learning} and distribution matching~\cite{liu2023generalized, yang2025topological}.


\subsection{Convergence rate with Marginal Constraints from an SDE}
The goal of this paper is to study the convergence rate of $\wht R^{\m T_m}$ to its limit as $m\to \infty$, when the marginal constraints $\mu_0, \dots, \mu_m$ in~\eqref{eqn: general_MSB} are derived from a stochastic differential equation (SDE).
More precisely, let $\{\rho_t\in\ms P(\m X): t\in[0, T]\}$ be a family of probability distributions satisfying the following assumption:
\begin{Assumption}\label{assump: SDE} 
$\rho_t$ is the marginal distribution at time $t$ of the stochastic process $Z_t$ defined through the SDE: 
\begin{align}\label{eqn: SDE}
    \dd Z_t = \nabla\Psi(t, Z_t)\,\dd t + \sqrt{\tau}\,\dd B_t,\quad Z_0\sim\rho_0\in\ms P_2^r(\m X),
\end{align}
where $\rho_0$ satisfies $\KL(\rho_0\,\|\,\lambda_{\m X}) < \infty$, $B_t$ is the standard Brownian motion on $\m X$, and $\Psi$ is a differentiable function on $\mathbb{R}_{\ge 0}\times \m X$, which we call a {\em potential function}.
Let $$R^\ast = \Law(\{Z_t: t\in[0, T]\})\in\ms P(\Omega)$$ be the distribution of the solution $Z$ of the above SDE. 
\end{Assumption}

By restricting the marginal constraints in~\eqref{eqn: general_MSB} to $\mu_j = \rho_{t_j}$ for $j = 0, 1, \dots, m$, the original multi-marginal Schrodinger bridge problem turns into
\begin{align}\label{eqn: SB_SDE}
\wht R^{\m T_m} = \argmin_{\substack{R\in\ms P(\Omega)\\ R_{t_j} = \rho_{t_j}, j=0,1,\dots, m}} \KL(R\,\|\,W^\tau).
\end{align}

For simplicty of argument, we focus on the case  $$\m X=\mathbb{T}^d.$$
Because the $d$-dimensional flat torus $\mathbb{T}^d$ is compact, the solution to 
\eqref{eqn: SB_SDE} exists, which is unique due to strict convexity of $R\mapsto \KL(R\, \|\, W^\tau)$.
Our main result establishes the explicit convergence rate of $\wht R^{\m T_m}$ towards $R^\ast$ under Assumption~\ref{assump: SDE} when the potential function $\Psi$ possess certain smoothness. To proceed, let $$\Delta_m := \sup_{j\in[m]}|t_j - t_{j-1}|.$$

\begin{theorem}\label{thm: main result}
Let $\m X = \mb T^d$ be the flat torus, and suppose Assumption \ref{assump: SDE} and $\Psi\in\m C^{2,4}([0, T]\times\m X)$.
Then, for the solution $\wht R^{\m T_m}$ to \eqref{eqn: SB_SDE} there exists constants $C_1, C_2 > 0$  depending only on $T$, $\tau$, and $\Psi$, such that
\begin{align*}
\KL(R^\ast\,\|\,\wht R^{\m T_m}) \leq \Big[\frac{3C_1}{2} + \sqrt{\frac{5C_1}{2}}C_2\Big]\Delta_m.
\end{align*}
When the time points are equally separated, we have $\KL(R^\ast\,\|\,\wht R^{\m T_m}) \lesssim m^{-1}$, where the constant of $\lesssim$ does not depend on $m$.
\end{theorem}


As we explain in \cref{sec:many-intervals}, \cref{thm: main result} is a consequence of the following result for the Schr\"odinger bridge problem ($m=1$):
\begin{theorem}[The $m=1$ case]\label{thm:epsilon-interval}
Make the same assumptions as in \cref{thm: main result}. 
For each $\epsilon >0$, consider  \begin{align}\label{eqn:single-SB}
 \wht R_{[0,\varepsilon]} =\argmin_{\substack{R\in\ms P(\Omega_{[0, \varepsilon]})\\R_0 = \rho_0, R_\epsilon =\rho_\epsilon}} \KL(R\,\|\,W^\tau_{[0,  \varepsilon]}). 
\end{align}
Then, there exists constants $C'_1, C'_2 > 0$  depending only on  $\tau$, and $\Psi$, such that
    \begin{align}\label{eqn:epsilon-interval}
\KL(R^\ast_{[0,\varepsilon]}\,\|\,\wht R_{[0,\varepsilon]}) 
&\le  \Big[\frac{3C'_1}{2} + \sqrt{\frac{5C'_1}{2}}C'_2\Big]\varepsilon^2.
\end{align}
\end{theorem}
The point of this result is the $\varepsilon^2$ estimate in \eqref{eqn:epsilon-interval}, which then gives the first order estimate (in $\Delta_m$) in \cref{thm: main result}. For a proof of Theorem~\ref{thm:epsilon-interval}, see \cref{sec:epsilon-interval},

\subsection{Related  literature}
In this section, we review existing results closely related to our main theorem. We begin with the following result, which demonstrates that under Assumption~\ref{assump: SDE}, the path measure $R^\ast$ can be exactly recovered from its marginal distributions by solving an entropy minimization problem with infinitely many marginal constraints.
\begin{proposition}[Theorem 2.1 by~\cite{lavenant2024toward}]\label{prop:varaational-R*}
Assume $\m X$ is a compact smooth Riemannian manifold without boundary and the potential function $\Psi$ is $\m C^2$-smooth. Under Assumption~\ref{assump: SDE}, we have
\begin{align*}
R^\ast = \argmin_{\substack{R\in\ms P(\Omega)\\ R_t = \rho_t, t\in[0,T]}} \KL(R\,\|\,W^\tau).
\end{align*}
\end{proposition}

In practice, solving an entropy minimization problem with infinitely many marginal constraints can be computationally challenging or even impractical. It is therefore natural to consider an approximation involving finitely many marginal constraints. The following result shows that the approximation via the multi-marginal Schrodinger bridge problem weakly converges to the ideal solution of the original entropy minimization problem.
\begin{proposition}[Theorem 3.15 by~\cite{mohamed2021schrodinger}] 
Let $\m X=\mb R^d$. Assume that $\m T_1\subset\m T_2 \subset \dots \subset \m T_m\subset\dots \subset[0, T]$ is an increasing sequence of finite subsets, and let $\m T = \cup_{n=1}^\infty \m T_n$. If there exists a unique solution
\begin{align*}
\wht R^{\m T} \coloneqq \argmin_{\substack{R\in\ms P(\Omega)\\ R_{t} = \rho_t,\forall\, t\in\m T}} \KL(R\,\|\,W^\tau)
\end{align*}
with $\KL(\wht R^{\m T}\,\|\, W^\tau) < \infty$,
then there exists a unique solution $\wht R^{\m T_m}$ of problem~\eqref{eqn: SB_SDE} for every $m$. Furthermore, $\wht R^{\m T_m}$ converges to $\wht R^{\m T}$ weakly and $\KL(\wht R^{\m T_m}\,\|\,W^\tau)\to\KL(\wht R^{\m T}\,\|\,W^\tau)$ as $m\to\infty$.
\end{proposition}
When $\m T$ is a dense subset of $[0, T]$ and Assumption~\ref{assump: SDE} holds, the above two propositions together (ignoring their different state space $\m X$) imply that the solution $\wht R^{\m T_m}$ weakly converges to $R^\ast$ as $m\to\infty$. 

We are interested in deriving an explicit convergence rate depending on the sequence of time points $\m T_m$. To the best of our knowledge, the only existing result in this direction is provided by~\cite{agarwal2024iterated}, which establishes an explicit convergence rate under the assumption that the potential function $\Psi$ is time-independent.

\begin{proposition}[Corollary of Theorem 1 by~\cite{agarwal2024iterated}]\label{prop:better-rate}
When Assumption~\ref{assump: SDE} holds with a time independent potential function $\Psi(x)$ satisfying additional mild conditions, it holds that
\begin{align*}
\lim_{m\to\infty}\Delta_m^{-1}\KL(R^\ast\,\|\,\wht R^{\m T_m}) = 0.
\end{align*}
\end{proposition}
In \cite[Theorem 1]{agarwal2024iterated}, they obtained a bound similar to our \cref{thm:epsilon-interval} for time-independent $\Psi$, but with a better rate $o(\epsilon^2)$, resulting in the order $o(\Delta_m)$ in their version of Proposition~\ref{prop:better-rate}. 
We note that the time-dependent case we consider in this paper is an important generalization; for example, it is relevant to statistical inference problems involving development trajectories in cell biology, e.g. \cite{lavenant2024toward, yao2025learning}, where the drift term for the developmental process is naturally time dependent; cells may follow different rules for development at different phases.

Our proof of the main results in Section~\ref{sec:main-proof} is also inspired by ~\citep{agarwal2024iterated}.
A crucial fact that we will use is the variational formulation of the Schr\"odinger bridge problem: 
First, notice that \eqref{eqn:single-SB} is equivalent to 
\begin{align}\label{eqn: scale_SB}
\wht R^\varepsilon = \argmin_{R_0 = \rho_0, R_1 = \rho_\varepsilon}\KL(R\,\|\,W^\varepsilon).
\end{align}
Then, when we focus on the marginal of $\wht R^\varepsilon$ we have
\begin{proposition}[Section 4.5 in~\cite{chen2021stochastic}]\label{prop:variational-SB}
    The problem \eqref{eqn: scale_SB} is then equivalent to the following:
\begin{align}\label{eqn:SB-dynamic}
&\inf_{(\mu_t, v_t)_{t\in[0, 1]}} \int_0^1\Big[\frac{1}{2}\|v_t\|^2_{L^2(\mu_t)}\, + \frac{\varepsilon^2}{8}\|\nabla\log\mu_t\|^2_{L^2(\mu_t)}\Big]\,\dd t,\\ \nonumber
&\mx{s.t.}\quad \partial_t\mu_t + \nabla\cdot(\mu_t v_t) = 0, \,\,\mu_0 = \rho_0,\,\,\mu_1 = \rho_\varepsilon,
\end{align}
where for the optimal solution $(\mu_t^\varepsilon, v_t^\varepsilon)$,  $\mu_t^\varepsilon$ is nothing but the marginal of $\wht R^\varepsilon$ at $t$.
\end{proposition}




\section{Proof of Main Results}\label{sec:main-proof}
The whole proof of our main result Theorem~\ref{thm: main result} consists of several steps. First, using the Markov property of the multi-marginal Schrodinger bridge $\wht R^{\m T_m}$, the associated KL divergence can be decomposed to the sum of KL divergence associated with Schrodinger bridges with only two marginal constraints (Section~\ref{sec:many-intervals}). 
Then, Section~\ref{sec:epsilon-interval} focuses on the two-marginal case. The KL divergence between the path measures of the SDE and the Schrodinger bridge can be controlled via their differences on the probability distributions of two marginals and the distribution of the paths connecting the marginals (see \textbf{Step 1} of the proof). The former can be controlled simply via Girsanov's Theorem (see \textbf{Step 2} of the proof), while the latter can be analyzed using the variational formulation of Schrodinger bridge (see \textbf{Step 3} of the proof). Lastly, combining all the pieces yields the desired estimate of the two-marginal case.



\subsection{Poof of \cref{thm: main result}  assuming the result of \cref{thm:epsilon-interval}}\label{sec:many-intervals}

First, applying Equation (10) in~\citep{leonard2012schrodinger} yields
\begin{align*}
\KL(R^\ast\,\|\,\wht R^{\m T_m})
&= \KL(R^\ast_{t_0, \cdots, t_m}\,\|\,\wht R^{\m T_m}_{t_0,\cdots, t_m}) + \mb E_{R^\ast_{t_0,\cdots, t_m}}\big[\KL\big(R^\ast(\cdot\,|\,Z_{t_0}, \cdots, Z_{t_m})\,\big\|\,\wht R^{\m T_m}(\cdot\,|\,Z_{t_0}, \cdots, Z_{t_m})\big)\big].
\end{align*}
Since $\wht R^{\m T_m}$ is the solution to~\eqref{eqn: SB_SDE}, we know that $\wht R^{\m T_m}$ is Markov (Proposition D.1 in~\cite{lavenant2024toward} or Lemma 3.4 in~\citep{benamou2019entropy}) and satisfies
\begin{align*}
\wht R^{\m T_m}_{t_0, \cdots, t_m}(\dd x_0, \cdots, \dd x_m) = \wht R^{\wht T_m}_{t_0, t_1}(\dd x_0, \dd x_1) \wht R^{T_m}_{t_2\,|\,t_1}(\dd x_2\,|\, x_1)\cdots \wht R^{T_m}_{t_m\,|\,t_{m-1}}(\dd x_m\,|\, x_{m-1}).
\end{align*}
Therefore, we have
\begin{align*}
\KL(R^\ast_{t_0, \cdots, t_m}\,\|\,\wht R^{\m T_m}_{t_0,\cdots, t_m}) = \sum_{j=1}^m \KL(R^\ast_{t_{j-1}, t_j}\,\|\,\wht R^{\m T_m}_{t_{j-1}, t_j}) - \sum_{j=1}^{m-1} \KL(R^\ast_{t_j}\,\|\,\wht R^{\m T_m}_{t_j})
= \sum_{j=1}^m \KL(R^\ast_{t_{j-1}, t_j}\,\|\,\wht R^{\m T_m}_{t_{j-1}, t_j}),
\end{align*}
where the last equation is due to the constraints $R_{t_j}^\ast = \wht R_{t_j}^{\m T_m} = \rho_{t_j}$ for all $j = 0, 1, \cdots, m$. Since both $R^\ast$ and $\wht R^{\m T_m}$ are Markov, it holds that
\begin{align*}
&\quad\,\KL\big(R^\ast(\cdot\,|\,Z_{t_0},\cdots, Z_{t_m})\,\big\|\,\wht R^{\m T_m}(\cdot\,|\,Z_{t_0},\cdots, Z_{t_m})\big)\\
&= \sum_{j=1}^m \KL\big(R^\ast_{[t_{j-1}, t_j]}(\cdot\,|\,Z_{t_{j-1}}, Z_{t_j})\,\big\|\,\wht R^{\m T_m}_{[t_{j-1}, t_j]}(\cdot\,|\,Z_{t_{j-1}}, Z_{t_j})\big).
\end{align*}
Therefore, we have
\begin{align*}
&\quad\,\mb E_{R^\ast_{t_0,\cdots, t_m}}\big[\KL\big(R^\ast(\cdot\,|\,Z_{t_0}, \cdots, Z_{t_m})\,\big\|\,\wht R^{\m T_m}(\cdot\,|\,Z_{t_0}, \cdots, Z_{t_m})\big)\big]\\
&\stackrel{}{=}\sum_{j=1}^m \mb E_{R^\ast_{t_0,\cdots, t_m}}\big[\KL\big(R^\ast_{[t_{j-1}, t_j]}(\cdot\,|\,Z_{t_{j-1}}, Z_{t_j})\,\big\|\,\wht R^{\m T_m}_{[t_{j-1}, t_j]}(\cdot\,|\,Z_{t_{j-1}}, Z_{t_j})\big)\big]\\
&= \sum_{j=1}^m \mb E_{R^\ast_{t_{j-1},t_j}}\big[\KL\big(R^\ast_{[t_{j-1}, t_j]}(\cdot\,|\,Z_{t_{j-1}}, Z_{t_j})\,\big\|\,\wht R^{\m T_m}_{[t_{j-1}, t_j]}(\cdot\,|\,Z_{t_{j-1}}, Z_{t_j})\big)\big]\\
&\stackrel{}{=} \sum_{j=1}^m \big[\KL\big(R^\ast_{[t_{j-1}, t_j]}\,\big\|\,\wht R^{\m T_m}_{[t_{j-1}, t_j]}\big) - \KL(R^\ast_{t_{j-1}, t_j}\,\|\,\wht R^{\m T_m}_{t_{j-1}, t_j})\big].
\end{align*}
Here, the last equality is due to Equation (10)~\citep{leonard2012schrodinger} again. Combining all pieces above yields
\begin{align}\label{eqn: MSB_piece}
\KL(R^\ast\,\|\,\wht R^{\m T_m}) = \sum_{j=1}^m \KL\big(R^\ast_{[t_{j-1}, t_j]}\,\|\, \wht R^{\m T_m}_{[t_{j-1}, t_j]}\big).
\end{align}
Moreover, we have from \cref{prop:varaational-R*} that
\begin{align*}
 \wht R^{\m T_m}_{[t_{j-1}, t_j]} = \argmin_{\substack{R\in\ms P(\Omega_{[t_{j-1}, t_j]})\\ R_{t_j} = \rho_{t_j}}} \KL(R\,\|\,W^\tau_{[t_{j-1}, t_j]}).
\end{align*}
Note that we can apply \cref{thm:epsilon-interval}  for each  terms in the summation of \eqref{eqn: MSB_piece} after a time shift. 
Therefore, we get
\begin{align*}
\KL(R^\ast\,\|\,\wht R^{\m T_m}) \
\leq \Big[\frac{3C'_1}{2} + \sqrt{\frac{5C'_1}{2}}C'_2\Big]\sum_{j=1}^m (t_j - t_{j-1})^2
\leq \Big[\frac{3C_1}{2} + \sqrt{\frac{5C_1}{2}}C_2\Big]\Delta_m
\end{align*}
 by setting $C_1=C_1'T$ and  $C_2=C_2'\sqrt{T}$.
 This completes the proof. \qed

\subsection{Proof of \cref{thm:epsilon-interval} (the $m=1$ case)}\label{sec:epsilon-interval}

We now prove \cref{thm:epsilon-interval} in four steps. From now on,
without loss of generality, we only need to consider the case where $$\tau = 1,$$ and we simply write $W^\tau$ as $W$. We also  abuse the notation and 
\begin{align*}
    \hbox{let  $W^{\rho_0}$ represent the Wiener measure with the initial distribution $W^{\rho_0}_0 = \rho_0$.} 
\end{align*}

\noindent\underline{{\bf Step 1:}}
We will prove that
\begin{align}\label{eqn: upp_bd}
\KL(R^\ast_{[0,\varepsilon]}\,\|\,\wht R_{[0,\varepsilon]}) \leq \KL(R^\ast_{[0,\varepsilon]}\,\|\,W^{\rho_0}_{[0,\varepsilon]}) + \mb E_{\wht R_{0,\varepsilon}}\Big[\log\frac{\dd W_{0, \varepsilon}^{\rho_0}}{\dd R^\ast_{0, \varepsilon}}\Big].
\end{align}

\begin{proof}
First, applying Equation (10) in~\citep{leonard2012schrodinger} yields
\begin{align*}
\KL(R^\ast_{[0,\varepsilon]}\,\|\,\wht R_{[0, \varepsilon]})
&= \KL(R^\ast_{0,\varepsilon}\,\|\,\wht R_{0,\varepsilon}) + \mb E_{R_{0,\varepsilon}^\ast}\Big[\KL\big(R_{[0,\varepsilon]}^\ast(\cdot\,|\,Z_0, Z_\varepsilon)\,\|\, \wht R_{[0,\varepsilon]}(\cdot\,|\,Z_0,Z_\varepsilon)\big)\Big],\\
\KL(R^\ast_{[0, \varepsilon]}\,\|\,W^{\rho_0}_{[0, \varepsilon]})
&= \KL(R_{0, \varepsilon}^\ast\,\|\,W_{0, \varepsilon}^{\rho_0}) + \mb E_{R_{0,\varepsilon}^\ast}\Big[\KL\big(R_{[0,\varepsilon]}^\ast(\cdot\,|\,Z_0, Z_\varepsilon)\,\|\, W^{\rho_0}_{[0,\varepsilon]}(\cdot\,|\,Z_0,Z_\varepsilon)\big)\Big].
\end{align*}
Combining with the fact that $$\wht R_{[0, \varepsilon]}(\cdot\,|\,Z_0, Z_\varepsilon) = W^{\rho_0}_{[0, \varepsilon]}(\cdot\,|\,Z_0, Z_\varepsilon)$$ (see e.g.~\citep[Proposition 2.3,][]{leonard2013survey}), we have
\begin{align*}
\KL(R^\ast_{[0,\varepsilon]}\,\|\,\wht R_{[0, \varepsilon]})
&= \KL(R_{[0,\varepsilon]}^\ast\,\|\,W^{\rho_0}_{[0,\varepsilon]}) +\KL(R^\ast_{0,\varepsilon}\,\|\,\wht R_{0,\varepsilon})   - \KL(R^\ast_{0,\varepsilon}\,\|\,W^{\rho_0}_{0,\varepsilon}).
\end{align*}
Note that we also have
\begin{align*}
\KL(R_{0,\varepsilon}^\ast\,\|\,\wht R_{0,\varepsilon}) - \KL(R_{0,\varepsilon}^\ast\,\|\,W_{0,\varepsilon}^{\rho_0})
&= \mb E_{R_{0,\varepsilon}^\ast}\Big[\log\frac{\dd R_{0,\varepsilon}^\ast}{\dd \wht R_{0,\varepsilon}} - \log\frac{\dd R_{0,\varepsilon}^\ast}{\dd W_{0,\varepsilon}^{\rho_0}}\Big] = \mb E_{R_{0,\varepsilon}^\ast}\Big[\log\frac{\dd W_{0,\varepsilon}^{\rho_0}}{\dd \wht R_{0,\varepsilon}}\Big].
\end{align*}
Let $p_\varepsilon(x, y)$ be the transition probability of $W_{0,\varepsilon}$. Then, we have
\begin{align*}
W_{0,\varepsilon}^{\rho_0}(x, y) 
&= \rho_0(x)p_\varepsilon(x, y).
\end{align*}
For the  Schrodinger bridge problem~\eqref{eqn:single-SB},
it is well-known (see e.g. \cite{leonard2012schrodinger, nutz2021introduction}) that 
there exist potential functions $\varphi$ and $\psi$ on $\m X$, such that
\begin{align*}
\wht R_{0,\varepsilon}(x, y) = e^{\frac{\varphi(x) + \psi(y)}{\varepsilon}} p_\varepsilon(x, y) \rho_0(x)\rho_\varepsilon(y).
\end{align*}
Therefore, we have
\begin{align*}
\mb \log\frac{\dd W_{0,\varepsilon}^{\rho_0}}{\dd\wht R_{0,\varepsilon}}(x, y) =  - \log\rho_\varepsilon(y) - \left(\frac{\varphi(x) + \psi(y)}{\varepsilon}\right),
\end{align*}
which is the sum of functions of $x$ and of $y$ without any term depending jointly on $(x,y)$.  Therefore, its  expected value with respect to ${R_{0,\varepsilon}^\ast}$ depends only on the marginals at $t=0$ and $t=\epsilon$ of ${R_{0,\varepsilon}^\ast}$.
As $\wht R_{0,\varepsilon}$ has the same marginals as ${R_{0,\varepsilon}^\ast}$, we have
\begin{align*}
\mb E_{R_{0,\varepsilon}^\ast}\Big[\log\frac{\dd W_{0,\varepsilon}^{\rho_0}}{\dd \wht R_{0,\varepsilon}}\Big]
& = \mb E_{\wht R_{0,\varepsilon}}\Big[\log\frac{\dd W_{0,\varepsilon}^{\rho_0}}{\dd \wht R_{0,\varepsilon}}\Big]\\
&= \mb E_{\wht R_{0,\varepsilon}}\Big[\log\frac{\dd W_{0, \varepsilon}^{\rho_0}}{\dd R^\ast_{0, \varepsilon}}\Big] -\KL(\wht R_{0,\varepsilon}\,\|\,R_{0, \varepsilon}^\ast)\\
& \leq \mb E_{\wht R_{0,\varepsilon}}\Big[\log\frac{\dd W_{0, \varepsilon}^{\rho_0}}{\dd R^\ast_{0, \varepsilon}}\Big]  \quad \hbox{(as $\KL \ge 0$)}
\end{align*}
This proves \eqref{eqn: upp_bd}.
   
\end{proof}

\bigskip

\noindent\underline{\bf Step 2}: 
We now prove that $\displaystyle \KL(R^\ast_{[0,\varepsilon]}\,\|\,W^{\rho_0}_{[0,\varepsilon]}) = \frac{1}{2}\int_0^\varepsilon \mb E\|\nabla\Phi(t, Z_t)\|^2\,\dd t$, where we recall that $\{Z_t: 0\leq t\leq T\}$ is the stochastic process in Assumption~\ref{assump: SDE}.

\begin{proof}
 
By Girsanov's Theorem~\citep[Theorem 1.12,][]{kutoyants2013statistical}
\begin{align*}
\log\frac{\dd R^\ast_{[0, \varepsilon]}}{\dd W_{[0, \varepsilon]}^{\rho_0}} (\omega)
= \int_0^\varepsilon\nabla\Psi(t, \omega_t)^\top\dd\omega_t - \frac{1}{2}\int_0^\varepsilon\|\nabla\Psi(t,\omega_t)\|^2\,\dd t.
\end{align*}

Therefore, we have
\begin{align*}
\KL(R^\ast_{[0,\varepsilon]}\,\|\,W^{\rho_0}_{[0,\varepsilon]})
&= \mb E_{R^\ast}\bigg[\int_0^\varepsilon\nabla\Psi(t, {\omega_t})^\top\dd \omega_t - \frac{1}{2}\int_0^\varepsilon\|\nabla\Psi(t,{\omega_t})\|^2\,\dd t\bigg]\\
& = \mb E\bigg[\int_0^\varepsilon\nabla\Psi(t, Z_t)^\top\dd Z_t - \frac{1}{2}\int_0^\varepsilon\|\nabla\Psi(t,Z_t)\|^2\,\dd t\bigg]\\
& =
\mb E\bigg[\int_0^\varepsilon\nabla\Psi(t, Z_t)^\top \left(\nabla\Psi(t, Z_t)\,\dd t + \dd B_t\right) - \frac{1}{2}\int_0^\varepsilon\|\nabla\Psi(t,Z_t)\|^2\,\dd t\bigg] \ \ \hbox{ (recall $\tau=1$)}\\
&
= \frac{1}{2}\int_0^\varepsilon \mb E\|\nabla\Phi(t, Z_t)\|^2\,\dd t \ \ \hbox{(as   $\mb E\left[ \int_0^\epsilon \nabla\Psi(t, Z_t)^\top  \dd B_t \right]=0$)}.
\end{align*}

\end{proof}
\begin{rem}\label{rmk:Girsanov-general}
Knowing that  $R^\ast$ and $\wht R$ have the same initial distribution, one may use the Girsanov theorem ~\citep[Theorem 1.12,][]{kutoyants2013statistical} to get
\begin{align}\label{eqn:Girsanov-ratio}
    \frac{\dd R^\ast_{[0, \varepsilon]}}{\dd \wht R_{[0, \varepsilon]}} (\omega)
= \exp \left( \int_0^\varepsilon\nabla (\Psi(t, \omega_t)-\phi(t, \omega_t))^\top\dd\omega_t - \frac{1}{2}\int_0^\varepsilon (\|\nabla(\Psi(t,\omega_t)\|^2- \|\nabla \phi(t, \omega_t))\|^2)\,\dd t\right),
\end{align}
where $\phi$ is the optimal potential to the solution $\wht R_{[0, \varepsilon]}$ of \eqref{eqn:epsilon-interval} that gives the drift term in the corresponding SDE.
We then follow the same argument as in Step 2,  to get:
\begin{align}\label{eqn:KL-R*-Rhat}
\KL ( R^\ast_{[0, \varepsilon]} \, \| \,  \wht R_{[0, \varepsilon]})    &=\mb E \left[ \frac{1}{2}\int_0^\varepsilon \|\nabla\Psi(t,Z_t)- \nabla \phi(t, Z_t)\|^2\,\dd t \right]
\end{align}
\noindent This would give us a way to get the $O(\varepsilon^2)$ of  $\KL ( R^\ast_{[0, \varepsilon]} \, \| \,  \wht R_{[0, \varepsilon]})$ directly, if we knew how to control the difference  $\nabla \Psi -\nabla \phi$. It is not clear to us how to proceed this way. Our method instead compares $R^*$ first with the Brownian motion.

    
\end{rem}

\bigskip

\noindent\underline{\bf Step 3:} We will prove that there exists constants $C_1, C_2$ depending only on $\Psi$ such that
\begin{align*}
\mb E_{\wht R_{0,\varepsilon}}\Big[\log\frac{\dd W_{0, \varepsilon}^{\rho_0}}{\dd R^\ast_{0, \varepsilon}}\Big] 
\leq -\frac{1}{2}\int_0^\varepsilon\mb E\|\nabla\Phi(t, Z_t)\|^2\,\dd t + C_1\varepsilon^2  + \Big[\frac{C_1}{2} + \sqrt{\frac{5C_1}{2}}C_2\Big]\varepsilon^2.
\end{align*}


\begin{proof}
To see this, first let $W_{[0, \varepsilon]}^x\in\ms P(\Omega)$ be the Wiener measure on time $[0, \varepsilon]$ with initial distribution $\delta_{x}$, and $(R^\ast_{[0,\varepsilon]})^x$ is defined similarly. Since $R^\ast_0 = \rho_0$, we have, from Girsanov's theorem,
\begin{align*}
\frac{\dd R_{0,\varepsilon}^\ast}{\dd W_{0,\varepsilon}^{\rho_0}}(x, y)
&= \frac{\dd R_{\varepsilon|0}^\ast}{\dd W_{\varepsilon|0}}(x, y)
= \mb E_{W^x}\Big[\frac{\dd (R_{[0,\varepsilon]}^\ast)^x}{\dd W_{[0,\varepsilon]}^x}\,\Big|\, \omega_\varepsilon = y\Big]\\
&= \mb E\bigg[\exp\bigg\{\int_0^\varepsilon\nabla\Psi(t, B_t^x)^\top\,\dd B_t^x - \frac{1}{2}\int_0^\varepsilon\|\nabla\Psi(t, B_t^x)\|^2\,\dd t\bigg\}\,\Big|\,B_\varepsilon^x = y\bigg],
\end{align*}
where $B_t^x$ is the Brownian motion starting at $B_0^x = x$. By Ito's formula, we have
\begin{align*}
\dd\Psi(t, B_t^x) = \Big[\partial_t\Psi(t, B_t^x) + \frac{1}{2}\Delta\Psi(t, B_t^x)\Big]\,\dd t + \nabla\Psi(t, B_t^x)^\top\dd B_t^x,
\end{align*}
which implies
\begin{align*}
\int_0^\varepsilon\nabla\Psi(t, B_t^x)^\top\dd B_t^x 
= \big[\Psi(\varepsilon, B_\varepsilon^x) - \Psi(0, x)\big] - \int_0^\varepsilon\left(\partial_t \Psi(t, B_t^x) +  \frac{1}{2}\Delta\Psi(t, B_t^x)\right)\,\dd t.
\end{align*}
So, we have
\begin{align*}
&\frac{\dd R_{0,\varepsilon}^\ast}{\dd W_{0,\varepsilon}^{\rho_0}}(x, y)
= \mb E\Big[\exp\Big\{\Psi(\varepsilon, B_\varepsilon^x)- \Psi(0, x) - \int_0^\varepsilon\m U(t, B_t^x)\,\dd t\Big\}\,\Big|\,B_\varepsilon^x = y\Big],
\end{align*}
where
\begin{align*}
 \m U(t, z) = \partial_t\Psi(t, z) + \frac{1}{2}\Delta\Psi(t, z) + \frac{1}{2}\|\nabla\Psi(t, z)\|^2.
\end{align*}
To sum up, we have  derived for $(X, Y) \sim \wht R_{0, \varepsilon}$,
\begin{align*}
& \mb E_{\wht R_{0,\varepsilon}}\Big[\log\frac{\dd W_{0, \varepsilon}^{\rho_0}}{\dd R^\ast_{0, \varepsilon}}\Big]
= -\mb E_{\wht R_{0,\varepsilon}}\bigg[\log\mb E\Big[\exp\Big\{\Psi(\varepsilon, B_\varepsilon^X)- \Psi(0, X) - \int_0^\varepsilon\m U(t, B_t^X)\,\dd t\Big\}\,\Big|\,B_\varepsilon^X = Y\Big]\bigg]\\
&= \mb E_{\wht R_{0,\varepsilon}}\bigg[\big[\Psi(0, X) - \Psi(\varepsilon, Y)\big] + \varepsilon\m U(0, X) - \log\mb E\Big[e^{-\int_0^\varepsilon \left( \m U(t, B_t^X) - \m U(0, X)\right)\,\dd t}\,\Big|\,B_\varepsilon^X = Y\Big]\bigg].
\end{align*}

\noindent \underline{\bf Step 3.1:} we will show that there exists a constant $C_1 = C_1(\Psi)$ only depending on $\Psi$, such that
\begin{align*}
\mb E_{\wht R_{0,\varepsilon}}\big[\Psi(0, X) - \Psi(\varepsilon, Y) + \varepsilon\m U(0, X)\big] 
\leq -\frac{1}{2}\int_0^\varepsilon\mb E\|\nabla\Psi(t, Z_t)  \|^2\,\dd t + C_1\varepsilon^2.
\end{align*}
In fact, we have
\begin{align*}
\mb E_{\wht R_{0,\varepsilon}}\big[\Psi(0, X) - \Psi(\varepsilon, Y)\big]
&= \mb E\Psi(0, Z_0) - \mb E\Psi(\varepsilon, Z_\varepsilon)
= -\mb E\int_0^\varepsilon \dd\Psi(t, Z_t),
\end{align*}
where the first equality follows from that $\wht R_{0,\varepsilon}$ and $R^*_{0,\varepsilon}$ have the same marginals.

By Ito's formula, we have
\begin{align*}
\dd \Psi(t, Z_t) 
&= \Big[\partial_t\Psi(t, Z_t) + \frac{1}{2}\Delta\Psi(t, Z_t)\Big]\,\dd t + \nabla\Psi(t, Z_t)^\top \dd Z_t\\
&= \Big[\partial_t\Psi(t, Z_t) + \frac{1}{2}\Delta\Psi(t, Z_t) + \|\nabla\Psi(t, Z_t)\|^2\Big]\,\dd t + \nabla\Psi(t, Z_t)^\top\dd B_t\\
&=\Big[\m U(t, Z_t) + \frac{1}{2}\|\nabla\Psi(t, Z_t)\|^2\Big]\,\dd t + \nabla\Psi(t, Z_t)^\top\dd B_t. 
\end{align*}
Therefore, we have
\begin{align*}
\mb E\int_0^\varepsilon\dd\Psi(t, Z_t)
&= \int_0^\varepsilon\m U(t, Z_t) + \frac{1}{2}\|\nabla\Psi(t, Z_t)\|^2\,\dd t.
\end{align*}
This implies
\begin{align*}
\mb E_{\wht R_{0,\varepsilon}}\big[\Psi(0, X) - \Psi(\varepsilon, Y) + \varepsilon\m U(0, X)\big] = -\frac{1}{2}\int_0^\varepsilon\mb E\|\nabla\Psi(t, Z_t)   \|^2\,\dd t - \mb E\int_0^\varepsilon\m U(t, Z_t) - \m U(0, Z_0)\,\dd t.
\end{align*}
Again, by Ito's formula, we have
\begin{align*}
\dd \m U(t, Z_t) 
&=  
\Big[\partial_t\m U(t, Z_t) + \frac{1}{2}\Delta\m U(t, Z_t)\Big]\,\dd t + \nabla\m U(t, Z_t)^\top\,\dd Z_t
\\
&  = \Big[\partial_t\m U(t, Z_t) + \frac{1}{2}\Delta\m U(t, Z_t) + \nabla\m U(t, Z_t)^\top \nabla \Psi (t, Z_t)\Big]\,\dd t + \nabla\m U(t, Z_t)^\top\,\dd B_t.
\end{align*}
So, we have
\begin{align*}
\mb E\int_0^\varepsilon\m U(t, Z_t) - \m U(0, Z_0)\,\dd t
&=\int_0^\varepsilon\!\!\int_0^t\mb E\Big[\partial_s\m U(s, Z_s) + \frac{1}{2}\Delta\m U(s, Z_s) + \nabla\m U(s, Z_s)^\top\nabla\Psi(s, Z_s)\Big]\,\dd s\dd t
\leq \frac{C_1(\Psi)}{2}\varepsilon^2,
\end{align*}
here, we take
\begin{align*}
C_1(\Psi) = \max_{t\in[0, 1], x\in\m X} \Big[\big|\partial_t\m U(t, x)\big| + \frac{1}{2}\big|\Delta\m U(t, x)\big| + |\nabla\m U(t, x)^\top\nabla\Psi(t, x)| + \|\nabla\m U(t, x)\|^2\Big].
\end{align*}  

\noindent \underline{\bf Step 3.2:} we will show that there exists a constant $C_2 = C_2(\Psi)$ only depending on $\Psi$, such that
\begin{align*}
\mb E_{\wht R_{0,\varepsilon}}\bigg[- \log\mb E\Big[e^{-\int_0^\varepsilon\m U(t, B_t^X) - \m U(0, X)\,\dd t}\,\Big|\,B_\varepsilon^X = Y\Big]\bigg] \leq \Big[\frac{C_1(\Psi)}{2} + \sqrt{\frac{5C_1(\Psi)}{2}}C_2(\Psi)\Big]\varepsilon^2.
\end{align*}

To address this, define the rescaled Brownian motion $B_t^\varepsilon = B_{{\varepsilon}t}$, and let $W^\varepsilon\in\ms P(\Omega)$ be the law of $B^\varepsilon$ with $B^\varepsilon_0 \sim\rho_0$. Then, we have
\begin{align*}
\mb E\Big[e^{-\int_0^\varepsilon\m U(t, B_t^X) - \m U(0, X)\,\dd t}\,\Big|\,B_\varepsilon^X = Y\Big]
= \mb E\Big[e^{-\varepsilon\int_0^1\m U(\varepsilon t, B_{t}^{\varepsilon, X}) - \m U(0, X)\,\dd t}\,\Big|\,B_1^{\varepsilon, X} = Y\Big].
\end{align*}
Also, consider the scaled Schrodinger bridge problem \eqref{eqn: scale_SB}:
\begin{align*}
\wht R^\varepsilon = \argmin_{R_0 = \rho_0, R_1 = \rho_\varepsilon}\KL(R\,\|\,W^\varepsilon).
\end{align*}
Then, we know  (see e.g.~\citep[Proposition 2.3,][]{leonard2013survey})  that  $$\wht R^\varepsilon(\dd\omega) = W^\varepsilon(\dd\omega\,|\,\omega_0, \omega_1) \wht R_{0,\varepsilon}(\dd\omega_0, \dd\omega_1).$$ This implies
\begin{align*}
&\quad\,\mb E_{\wht R_{0,\varepsilon}}\bigg[-\log\mb E\Big[e^{-\varepsilon\int_0^1\m U(\varepsilon t, B_{t}^{\varepsilon, X}) - \m U(0, X)\,\dd t}\,\Big|\,B_1^{\varepsilon, X} = Y\Big]\bigg]\\
&\leq \mb E_{\wht R_{0,\varepsilon}}\bigg[\mb E\Big[{\varepsilon\int_0^1\m U(\varepsilon t, B_{t}^{\varepsilon, X}) - \m U(0, X)\,\dd t}\,\Big|\,B_1^{\varepsilon, X} = Y\Big]\bigg] \quad \hbox{(from Jensen's inequality)}\\
&= \varepsilon\mb E_{\wht R^\varepsilon}\int_0^1\m U(\varepsilon t, X_t) - \m U(0, X_0)\,\dd t.
\end{align*}
Now, let us control the integrand. Assume $(\mu_t^\varepsilon, v_t^\varepsilon)$ solve the  variational formulation \eqref{eqn:SB-dynamic}  in Proposition~\ref{prop:variational-SB}.
Then, we know $\wht R^\varepsilon_t = \mu_t^\varepsilon$, and the pair 
$ (\mu_t, v_t)= (\rho_{\varepsilon t}, \varepsilon\nabla\Psi(\varepsilon t, \cdot) - \frac{\varepsilon}{2}\nabla\log\rho_{\varepsilon t})$  is admissible for \eqref{eqn:SB-dynamic}. Thus,
\begin{align}\label{eqn:variational-with-Psi}
& \int_0^1\Big[\frac{1}{2}\|v_t^\varepsilon\|^2_{L^2(\mu_t^\varepsilon)}\, + \frac{\varepsilon^2}{8}\|\nabla\log\mu_t^\varepsilon\|^2_{L^2(\mu_t^\varepsilon)}\Big]\,\dd t \\ \nonumber
& \leq \int_0^1\Big[\frac{\varepsilon^2}{2}\Big\|\nabla\Psi(\varepsilon t,\cdot) - \frac{1}{2}\nabla\log\rho_{\varepsilon t}\Big\|^2_{L^2(\rho_{\varepsilon t})}\, + \frac{\varepsilon^2}{8}\|\nabla\log\rho_{\varepsilon t}\|^2_{L^2(\rho_{\varepsilon t})}\Big]\,\dd t.
\end{align}

Now, we have
\begin{align*}
&\mb E_{\wht R^\varepsilon}\m U(\varepsilon t, X_t) - \m U(0, X_0)
\\ &= \int_0^t\frac{\dd}{\dd s}\mb E_{\mu_s^\varepsilon}\m U(\varepsilon s, X_s)\,\dd s\\
&= \int_0^t\frac{\dd}{\dd s}\int\m U(\varepsilon s, x)\mu_s^\varepsilon(x)\,\dd x\dd s\\
&= \varepsilon\int_0^t\mb E_{\mu_s^\varepsilon}[\partial_t\m U(\varepsilon s, X_s)]\,\dd s - \int_0^t\!\!\int\m U(\varepsilon s, x)\nabla\cdot(\mu_s^\varepsilon v_s^\varepsilon)\,\dd x\dd s\\
&= \varepsilon\int_0^t\mb E_{\mu_s^\varepsilon}[\partial_t\m U(\varepsilon s, X_s)]\,\dd s + \int_0^t\mb E_{\mu_s^\varepsilon}\big[\nabla\m U(\varepsilon s, X_s)^\top v_s^\varepsilon(X_s)\big]\,\dd s \qquad \hbox{(by integration by parts)}.
\end{align*}
Note that the first term can be bounded as
\begin{align*}
\varepsilon\int_0^t\mb E_{\mu_s^\varepsilon}[\partial_t\m U(\varepsilon s, X_s)]\,\dd s \leq \varepsilon t C_1(\Psi).
\end{align*}
For the second term, we use \eqref{eqn:variational-with-Psi} and get
\begin{align*}
&\quad\,\int_0^t\mb E_{\mu_s^\varepsilon}\big[\nabla\m U(\varepsilon s, X_s)^\top v_s^\varepsilon(X_s)\big]\,\dd s
\leq \int_0^1 \mb E_{\mu_s^\varepsilon}\|\nabla\m U(\varepsilon s, X_s)\|\cdot\|v_s^\varepsilon(X_s)\|\,\dd s\\
&\leq \sqrt{\int_0^1\mb E_{\mu_s^\varepsilon}\|\nabla \m U(\varepsilon s, X_s)\|^2\,\dd s}\cdot\sqrt{\int_0^1\mb E_{\mu_s^\varepsilon}\|v_s^\varepsilon(X_s)\|^2\,\dd s}\\
&\leq \sqrt{C_1(\Psi)}\cdot \varepsilon\sqrt{\int_0^1\Big\|\nabla\Psi(\varepsilon t,\cdot) - \frac{1}{2}\nabla\log\rho_{\varepsilon t}\Big\|^2_{L^2(\rho_{\varepsilon t})}\, + \frac{1}{4}\|\nabla\log\rho_{\varepsilon t}\|^2_{L^2(\rho_{\varepsilon t})} - \frac{1}{4}\|\nabla\log\mu_t^\varepsilon\|^2_{L^2(\mu_t^\varepsilon)}\,\dd t}\\
&\leq \varepsilon\sqrt{\frac{5C_1(\Psi)}{2}}C_2(\Psi),
\end{align*}
where
\begin{align*}
C_2(\Psi) = \sup_{t\in[0, 1], x\in\m X} \max\big\{|\nabla\Psi(t, x)|, |\nabla\log\rho(t, x)|\big\} < \infty
\end{align*}
due to the compactness of $[0, 1]\times \m X$. Notice that $C_2(\Psi)$ depends only on $\Psi$ as $\rho_t$ is the marginal distribution of $Z_t$ (which follows \eqref{eqn: SDE}. Thus, we get
\begin{align*}
\varepsilon\mb E_{\wht R^\varepsilon}\left[\int_0^1\left(\m U(\varepsilon t, X_t) - \m U(0, X_0)\right)\,\dd t \right]
& \leq \varepsilon\int_0^1 \left(\varepsilon t C_1(\Psi) + \varepsilon\sqrt{\frac{5C_1(\Psi)}{2}}C_2(\Psi)\right)\,\dd t\\
& = \Big[\frac{C_1(\Psi)}{2} + \sqrt{\frac{5C_1(\Psi)}{2}}C_2(\Psi)\Big]\varepsilon^2.
\end{align*}
\end{proof}

\bigskip

\noindent\underline{\bf Step 4: Conclusion.} Combining all above pieces yields
\begin{align*}
\KL(R^\ast_{[0,\varepsilon]}\,\|\,\wht R_{[0,\varepsilon]}) 
&\leq  \frac{1}{2}\int_0^\varepsilon \mb E\|\nabla\Phi(t, Z_t)\|^2\,\dd t + \bigg[-\frac{1}{2}\int_0^\varepsilon\mb E\|\nabla\Phi(t, Z_t)\|^2\,\dd t + C_1\varepsilon^2 + \Big[\frac{C_1}{2} + \sqrt{\frac{5C_1}{2}}C_2\Big]\varepsilon^2\bigg]\\
&= \Big[\frac{3C_1}{2} + \sqrt{\frac{5C_1}{2}}C_2\Big]\varepsilon^2.
\end{align*} \qed


\bibliography{ref}

\bibliographystyle{plain}

\end{document}